\documentclass[11pt]{article}
\usepackage{amsfonts}
\usepackage{color}
\newtheorem{theo}{Theorem}[section]
\newtheorem{lem}[theo]{Lemma}
\newtheorem{cor}[theo]{Corollary}

\newcommand{\mysection}[1]{\section{#1} \setcounter{equation}{0}}
\newcommand{\proof}{{\sc Proof.} \quad}
\newcommand{\proofc}{{\sc Proof} \ }
\newcommand{\be}{\begin{equation} \label}
\newcommand{\ee}{\end{equation}}
\newcommand{\bea}{\begin{eqnarray}\label}
\newcommand{\eea}{\end{eqnarray}}
\newcommand{\bas}{\begin{eqnarray*}}
\newcommand{\eas}{\end{eqnarray*}}
\newcommand{\bit}{\begin{itemize}}
\newcommand{\eit}{\end{itemize}}
\newcommand{\qed}{\hfill$\Box$ \vskip.2cm}
\newcommand{\nn}{\nonumber}
\newcommand{\R}{\mathbb{R}}
\newcommand{\N}{\mathbb{N}}
\newcommand{\pO}{\partial\Omega}
\newcommand{\dN}{\partial_\nu}
\newcommand{\eps}{\varepsilon}

\newcommand{\abs}{\\[5pt]}

\newcommand{\tm}{T_{max}}
\newcommand{\io}{\int_\Omega}

\newcommand{\hu}{{\hat{u}}}
\newcommand{\tu}{{\tilde{u}}}
\begin{document}
\title{Boundedness in a quasilinear parabolic-parabolic Keller-Segel system with
subcritical sensitivity}
\author{
Youshan Tao\\       
{\small Department of Applied Mathematics, Dong Hua University,}\\
{\small Shanghai 200051, P.R.~China}\\
{\small \tt taoys@dhu.edu.cn}
\and
Michael Winkler\\   
{\small Institut f\"ur Mathematik, Universit\"at Paderborn,}\\
{\small 33098 Paderborn, Germany}\\
{\small \tt michael.winkler@math.uni-paderborn.de} }
\date{}
\maketitle
\begin{abstract}
\noindent We consider the quasilinear parabolic-parabolic
Keller-Segel system
$$
    \left\{ \begin{array}{l}
    u_t=\nabla \cdot (D(u)\nabla u) - \nabla \cdot (S(u)\nabla v),
    \qquad x\in\Omega, \ t>0, \\[1mm]
    v_t=\Delta v -v + u,
    \hspace{3.55cm} x\in\Omega, \ t>0,
    \end{array} \right.
$$
under homogeneous Neumann boundary conditions in a smooth bounded
domain $\Omega\subset\R^n$ with $n\ge 2$.\\
\noindent It is proved that if $\frac{S(u)}{D(u)}\le cu^{\alpha}$
with $\alpha<\frac{2}{n}$ and some constant $c>0$ for all $u>1$ and
some further technical conditions are fulfilled, then the classical
solutions to the above system are uniformly-in-time bounded.
This boundedness result is optimal
according to a recent result by the second author
({\em Math. Meth. Appl. Sci.} {\bf 33} (2010), 12-24),
which says that if
$\frac{S(u)}{D(u)} \ge cu^\alpha$ for $u>1$ with $c>0$ and some
$\alpha>\frac{2}{n}$, then for each mass $M>0$ there exist blow-up solutions
with mass $\io u_0=M$.\\
\noindent In addition, this paper also proves a general boundedness
result for quasilinear non-uniformly parabolic equations by
modifying the iterative technique of Moser-Alikakos (Alikakos, {\em Comm.
PDE} {\bf 4} (1979), 827-868).\abs
\end{abstract}
%
%
%
%
%
%
\section*{Introduction}
This work is concerned with the initial-boundary value problem
\be{0}
    \left\{ \begin{array}{l}
    u_t=\nabla \cdot (D(u)\nabla u) - \nabla \cdot (S(u)\nabla v),
    \qquad x\in\Omega, \ t>0, \\[1mm]
    v_t=\Delta v -v + u,
    \hspace{3.9cm} x\in\Omega, \ t>0, \\[1mm]
    \dN u=\dN v =0, \hspace{4.2cm} x\in\partial\Omega, \ t>0, \\[1mm]
    u(x,0)=u_0(x),\quad  v(x,0)=v_0(x), \hspace{1.1cm} x\in\Omega,
    \end{array} \right.
\ee for the unknown $u=u(x,t), v=v(x,t)$, where $\Omega$ is a
bounded convex domain in $\R^n$ with smooth boundary, and $n\ge 2$.
The initial distributions $u_0$ and
$v_0$ are assumed to be nonnegative functions subject to the
inclusions $u_0 \in C^0(\bar\Omega)$ and $v_0 \in C^1(\bar\Omega)$,
respectively.\abs
Chemotaxis, the biassed movement of cells (or organisms) in response
to chemical gradients, plays an important role coordinating cell
migration in many biological phenomena (cf.~then review article
\cite{hillen_painter2009}). In (\ref{0}), $u$ denotes the cell
density and $v$ describes the concentration of the chemical signal
secreted by cells. In addition to (self-)diffusion, cells move
towards higher signal concentration, whereas the chemical signal undergoes
random diffusion and decay.
An important variant of the quasilinear
chemotaxis model (\ref{0}) was initially proposed by Painter and
Hillen \cite{painter_hillen}. Their approach assumes the presence
of a so-called {\em volume-filling effect}:
The movement of cells is inhibited near points where the cells are densely packed.
Painter and Hillen \cite{painter_hillen} derived their model
via random walk approach and they found a functional link between
the self-diffusivity $D(u)$ and the chemotactic sensitivity $S(u)$
that, in a non-dimensionalized version, takes the form
 \be{1}
    D(u)=Q(u)-uQ'(u),\qquad S(u)=uQ(u)
 \ee
where $Q(u)$ denotes the density-dependent probability for a cell to
find space somewhere in its neighboring location. Since this
probability is basically unknown, different choices for $Q$ are
conceivable.\abs
If $Q(u)\equiv 1$ we arrive at the classical Keller-Segel model
(\cite{keller_segel}),
 \be{2}
        \left\{ \begin{array}{ll}
        u_t=\Delta u - \nabla \cdot ( u\nabla v), \qquad &
    x\in\Omega, \ t>0, \\[1mm]
    v_t=\Delta v -v + u, &
    x\in\Omega, \ t>0, \\[1mm]
        \dN u=\dN v =0,  &
    x\in\partial\Omega, \ t>0, \\[1mm]
        u(x,0)=u_0(x), \quad v(x,0)=v_0(x),  \quad & x\in\Omega,
    \end{array} \right.
 \ee
 which has been investigated quite thoroughly during the past three decades.
 In view of the biologically meaningful question whether or not cell populations
 spontaneously form aggregates, most
 mathematical studies focused on whether solutions remain bounded or blow up.
 If $n=1$, then all solutions
 of (\ref{2}) are global in time and bounded (\cite{osaki_yagi}); if $n=2$ and $\io u_0<4\pi$, then the
 solution will be global and bounded (\cite{nagai_senba_yoshida}); if $n\ge 3$ and, for any $\delta>0$, the
 quantities $\| u_0\|_{L^{n/2+\delta}(\Omega)}$
 and $\| \nabla v_0\|_{L^{n+\delta}(\Omega)} $ are small, then the solution is
 global and bounded (\cite{winkler_jde}). On the other hand, if $n=2$ then for almost every $M>4\pi$ there exist smooth
 initial data $(u_0, v_0)$ with $\io u_0=M$ such that corresponding
 solution of (\ref{2}) blows up either in finite or infinite time
 provided $\Omega$ is simply connected (\cite{horstmann_wang}); in the particular framework
 of radially symmetric solutions in a planar disk, solutions may even blow
 up in finite time (\cite{herrero_velazuez}); if $n\ge 3$ and $\Omega$
 is a ball, then for all $M>0$ there exist initial data with $\io
 u_0=M$ such that the solution will become unbounded either in finite or infinite time
 (\cite{winkler_jde}).  \abs
In \cite{hillen_painter2001} the authors analyze (\ref{0}) upon the particular choices
$D(u)\equiv 1$ and $S(u)=u(1-u)_+$. This corresponds to the case of the compactly supported
probability $Q(u)=(1-u)_+$ in the volume-filling model, in particular meaning that the
chemotactic movement is entirely stopped when the cell density reaches the critical level $u=1$.
The resulting system admits global
bounded solutions only (\cite{hillen_painter2001}). Furthermore,
Wrzosek (\cite{wrzosek}, \cite{wrzosek2004}) studied the dynamical properties such as
instability of constant steady states or the existence of
attractors.\abs
The focus of this paper is to provide some further step towards a more detailed
understanding the interaction of the competing nonlinear mechanisms of
(self-)diffusion and cross-diffusion in (\ref{0}), allowing for rather general choices of $D(u)$ and $S(u)$.
Here we concentrate on the particular phenomenon of blow-up, and observe that in this respect,
previous results indicate that the asymptotic behavior of the ratio
$\frac{S(u)}{D(u)}$
for large values of $u$ seems to be decisive: Namely, in \cite{win_mmas} it has been shown that
\be{blowup}
    \begin{array}{l}
    \mbox{if } \quad
    \frac{S(u)}{D(u)} \ge c u^{\frac{2}{n}+\eps} \quad \mbox{for all } u>1 \mbox{ and some $c>0$ and $\eps>0$,} \\
    \hspace*{34mm}
    \mbox{then there exist smooth solutions of (\ref{0}) which blow up}
    \end{array}
\ee
either in finite or infinite time, provided that $\Omega$ is a ball.
%
However, to the best of our knowledge the existing literature leaves open the question
in how far this growth condition is critical in respect of blow-up.\\
It is the purpose of the present work to close this gap,
and correspondigly we shall suppose throughout that $D$ and $S$, besides
\be{reg}
    D \in C^2([0,\infty)) \qquad \mbox{and} \qquad S \in C^2([0,\infty)) \quad \mbox{with } S(0)=0,
\ee
are such that their ratio satisfies the growth condition
\be{alpha}
    \frac{S(u)}{D(u)} \le K (u+1)^\alpha \qquad \mbox{for all } u \ge 0
\ee
with some $K>0$ and $\alpha>0$.
Moreover, our approach will require the further technical assumptions that
\be{m}
    D(u) \ge K_0 (u+1)^{m-1} \quad \mbox{for all } u \ge 0
\ee
and
\be{M}
    D(u) \le K_1 (u+1)^{M-1} \quad \mbox{for all } u\ge 0
\ee
are valid with some constants $m\in\R, M\in \R, K_0>0$ and $K_1>0$.\abs
Under these hypotheses, our main result reads as follows.
\begin{theo}\label{theo_main}
  Suppose that $\Omega\subset \R^n$, $n\ge 2$, is a bounded convex domain with smooth boundary.
  Assume that $D$ and $S$ satisfy (\ref{reg}), (\ref{alpha}), (\ref{m}) and (\ref{M}) with some $m\in\R, M\in\R$
  and positive constants $K, K_0, K_1$ and
  \bas
    \alpha<\frac{2}{n}.
  \eas
  Then for any nonnegative $u_0 \in C^0(\bar\Omega)$ and $v_0 \in C^1(\bar\Omega)$,
  there exists a couple $(u,v)$ of nonnegative bounded functions belonging to
  $C^0(\bar\Omega \times [0,\infty)) \cap C^{2,1}(\bar\Omega \times (0,\infty))$ which solve (\ref{0}) classically.
\end{theo}
In conjunction with (\ref{blowup}), this provides an essentially complete picture on the dichotomy {\em boundedness
vs.~blow-up} in (\ref{0}), provided that the
(self-)diffusivity $D(u)$ has an asymptotically algebraic behavior.
It is an interesting open question that unfortunately has to be left open here whether the above boundedness statement
is also valid when $D(u)$ is allowed to grow or decay exponentially, for instance.\abs
Let us mention some further previous contributions in this direction.
In the particular case $D(u) \equiv 1$, the criticality of $\frac{S(u)}{D(u)} \simeq u^\frac{2}{n}$ was already revealed in
\cite{horstmann_winkler}, where global boundedness of solutions was shown when
$S(u) \le cu^{\frac{2}{n}-\eps}$ for all $u>1$ and some $c>0$ and $\eps>0$, and where some radial
blow-up solutions were constructed if $S(u) \ge cu^{\frac{2}{n}+\eps}$ for $u>1$ with $c>0$ and $\eps>0$, and if some further
technical restrictions hold.\\
As to the special case when $S(u)=u$,
Kowalczyk and Szyma\'{n}ska (\cite{kowalczyk_szymanska}) proved that solutions remain bounded under the condition
that $D(u) \ge c u^{2-\frac{4}{n}+\eps}$ for all $u>0$ with some $c>0$ and $\eps>0$.
In view of the above results, this is optimal for non-degenerate diffusion (with $D>0$ on $[0,\infty)$)
if and only if $n=2$.
For the same choice of $S(u)$ and $D>0$ on $[0,\infty)$,
Senba and Suzuki (\cite{senba_suzuki}) reached the critical exponent by showing boundedness
under the hypothesis that $D(u) \ge cu^{\frac{n-2}{n}+\eps}$ be valid for $u>1$ with some $c>0$ and $\eps>0$.
For more general $D(u)$ and $S(u)$ satisfying some technical assumptions, Cie\'slak (\cite{cieslak2008}) asserted
boundedness of solutions when either $n=2$ and $\frac{S(u)}{D(u)} \le cu^{\frac{1}{2}-\eps}$, or $n=3$ and
$\frac{S(u)}{D(u)} \le cu^{-\eps}$ for all $u>1$ and some $c>0$ and $\eps>0$ (cf.~also \cite{cieslak2007} for related
results).\abs
When the diffusion of the chemical signal is considered to occur much faster than
that of cells, by the approach of quasi-steady-state approximation
(cf.~\cite{jager_lackaus} or \cite{perthame}), the
parabolic-parabolic chemotaxis model (\ref{0}) can be reduced to simplified parabolic-elliptic models
where the second PDE in (\ref{0}) is replaced with either $0=\Delta v - v + u$, or with
$0=\Delta v - M + u$, where $M:=\io u_0$ denotes the total mass of cells.
For the former model,
if $n=2$, $S(u)=u$ and $D(u)\ge c(1+u)^{1+\eps}$ with $c>0$ and $\eps>0$, boundedness of solutions was proved in
\cite{kowalczyk}, and
the same conclusion was found in \cite{cieslak_moralesrodrigo} for more general $D(u)$ and $S(u)$
with the property that for some $c>0$ and $\eps>0$ we have
$\frac{S(u)}{D(u)} \le cu^{-\eps}$ when $n=2$, and $\frac{S(u)}{D(u)} \le cu^{-1-\eps}$ when $n=3$.\abs
As to the latter simplification, the knowledge appears to be rather complete and consistent with the results for the
parabolic-parabolic case if
$D(u) \simeq u^{-\gamma}$ and $S(u) \simeq u^\alpha$ for large $u$ with some $\gamma \ge 0$ and $\alpha\in\R$:
Solutions remain bounded if $\alpha+\gamma<\frac{2}{n}$, whereas blow-up may occur if $\alpha+\gamma>\frac{2}{n}$
(\cite{djie_winkler}, cf.~also \cite{cieslak_winkler} for a precedent addresing the special case $S(u)=u$).
Moreover, if $S(u)=u$, then even the critical case $D(u) \simeq u^\frac{n-2}{n}$ can be analyzed, and Cie\'slak and
Lauren\c{c}ot have shown it to belong to the blow-up regime (\cite{cieslak_laurencot}).
Refined conditions ensuring boundedness in two-dimensional parabolic-elliptic Keller-Segel models can be found in
\cite{calvez_carrillo}. For results in the whole space $\R^n$ with $D(u)$ and $S(u)$ being exact powers of $u$ (thus involving
porous medium-type or fast diffusion), we refer to \cite{sugiyama} and the references therein.\abs
The proof of our main results will be based on
a priori estimates in spatial Lebesgue spaces for $u$ and $\nabla v$.
Due to the careful adjustment of some parameters (cf.~Section \ref{sect_parameters}), our technique of deriving
integral bounds (see Section \ref{mainproof}) does not need
any iterative argument to establish bounds for $u(\cdot,t)$ in $L^p(\Omega)$ for any finite $p$,
as required in
some previous approaches (cf.~\cite{horstmann_winkler}, for instance).
Only in a final step an iteration is needed in order to turn this into a bound in $L^\infty(\Omega)$
by means of a Moser-Alikakos-type procedure (cf.~the appendix).\abs
%
%
%
%
%
%
%
%
\mysection{Local existence}\label{sect_prelim}
The following statement concerning local existence of classical solution can be proved by well-established methods
involving standard parabolic regularity theory and an appropriate fixed point framework
(for details see \cite{horstmann_winkler}, \cite{wrzosek2004} or also \cite{cieslak2007}, for instance).
\begin{lem}\label{lem_loc_exist}
  Let $D$ and $S$ satisfy (\ref{reg}), (\ref{m}) and (\ref{M}) with some $m\in\R,M\in\R, K_0>0$ and $K_1>0$,
  and assume that $u_0 \in C^0(\bar\Omega)$ and $v_0 \in C^1(\bar\Omega)$ are nonnegative.
  Then there exist $\tm \in (0,\infty]$ and a
  pair $(u,v)$ of functions from $C^0(\bar\Omega \times [0,\tm)) \cap C^{2,1}(\bar\Omega \times (0,\tm))$
  solving (\ref{0}) classically in $\Omega \times (0,\tm)$. These functions satisfy the inequalities
  \bas
    u \ge 0 \qquad \mbox{and} \qquad v \ge 0
    \qquad \mbox{in } \Omega \times (0,\tm),
  \eas
  and moreover
  \be{extend}
    \mbox{either $\tm=\infty$, \quad or }
    \limsup_{t\nearrow \tm} \Big( \|u(t)\|_{L^\infty(\Omega)} + \|v(t)\|_{L^\infty(\Omega)} \Big) = \infty.
  \ee
\end{lem}
The following properties of solutions of (\ref{0}) are well-known.
\begin{lem}\label{lem_basic}
  i) \ The first component $u$ of the solution of (\ref{0}) satisfies the mass conservation property
  \be{mass}
    \|u(t)\|_{L^1(\Omega)} = \|u_0\|_{L^1(\Omega)} \qquad \mbox{for all } t \in (0,\tm).
  \ee
  ii) \ For all $s\in [1,\frac{n}{n-1})$ there exists $c>0$ such that
  \be{gradv}
    \|v(t)\|_{W^{1,s}(\Omega)} \le c \qquad \mbox{for all } t\in (0,\tm)
  \ee
  holds.
\end{lem}
\proof
  Integrating with respect to $x\in\Omega$, we see that $\frac{d}{dt} \io u \equiv 0$, and that
  $\frac{d}{dt} \io v=-\io v + \io u$ for $t \in (0,\tm)$. This yields (\ref{mass}) and moreover shows that
  $v$ is bounded in $L^\infty((0,\tm);L^1(\Omega))$. Now this implies (\ref{gradv}) upon
  a standard regularity argument involving the variation-of-constants formula for $v$ and $L^p-L^q$ estimates
  for the heat semigroup (see \cite[Lemma 4.1]{horstmann_winkler}, for instance).
\qed
\mysection{Adjusting some parameters}\label{sect_parameters}
We now make sure that when the parameter $\alpha$ in (\ref{alpha}) indeed satisfies $\alpha<\frac{2}{n}$, we can choose
certain parameters, to be used in Lemma \ref{lem33} below, appropriately.
\begin{lem}\label{lem45}
  Let $n \ge 2, m\in\R, \alpha \in (0,\frac{2}{n}), \bar p \ge 1$ and $\bar q \ge 2$.
  Then there exist numbers $p\ge \bar p, q\ge \bar q, s\in [1,\frac{n}{n-1}), \theta>1$ and $\mu>1$ such that
  \begin{eqnarray}
    & & p>\max \Big\{ 4-m, \frac{n(1-m)}{2} \Big\}, \label{33.7} \\
    & & \frac{n-2}{n} \cdot \frac{m+p+2\alpha-3}{m+p-1} < \frac{1}{\theta}, \label{33.5} \\
    & & \frac{n-2}{n} \cdot \frac{2}{m+p-1} < \frac{1}{\mu}, \label{33.6} \\
    & & \frac{1}{\theta} < 1- \frac{n-2}{n} \cdot \frac{1}{q} \qquad \mbox{and} \label{33.13} \\
    & & \frac{1}{\mu} < \frac{2}{n} + \frac{n-2}{n} \cdot \frac{1}{q}, \label{33.14}
  \end{eqnarray}
  and such that moreover
  \be{33.1551}
    \frac{m+p+2\alpha-3-\frac{1}{\theta}}{1-\frac{n}{2}+\frac{n(m+p-1)}{2}}
    +\frac{\frac{2}{s}-1+\frac{1}{\theta}}{1-\frac{n}{2}+\frac{nq}{s}} < \frac{2}{n}
  \ee
  as well as
  \be{33.1552}
    \frac{2-\frac{1}{\mu}}{1-\frac{n}{2}+\frac{n(m+p-1)}{2}}
    + \frac{\frac{2(q-1)}{s}-1+\frac{1}{\mu}} {1-\frac{n}{2}+\frac{nq}{s}} < \frac{2}{n}
  \ee
  hold.
\end{lem}
\proof
  Let us first fix numbers $\theta>1$ and $\mu>1$ such that
  \be{33.3}
    (n-2)\theta < n
  \ee
  and
  \be{33.4}
    \mu>\frac{n}{2},
  \ee
  and let
  \bas
    q_0(p):=\frac{n(m+p-1)}{2(n-1)} \qquad \mbox{for } p\ge 1.
  \eas
  Then we can easily find some large $p\ge \bar p$ fulfilling
  \be{33.44}
    q_0(p) > \bar q,
  \ee
  and such that (\ref{33.7}), (\ref{33.5}) and (\ref{33.6}) hold as well as
  \be{33.8}
    \frac{1}{\theta} < 1 - \frac{n-2}{n} \cdot \frac{1}{q_0(p)}
  \ee
  and
  \be{33.9}
    \frac{1}{\mu} < \frac{2}{n} + \frac{n-2}{n} \cdot \frac{1}{q_0(p)}.
  \ee
  Here we note that (\ref{33.3}) asserts that (\ref{33.5}) is true for all sufficiently large $p$, whereas the fact that
  $q_0(p) \to + \infty$ as $p\to\infty$ along with the inequality $\theta>1$ and (\ref{33.4}) guarantees the validity
  of (\ref{33.8}) and (\ref{33.9}) for appropriately large $p$.\\
  We next let
  \be{33.99}
    f(q,s):=\frac{m+p+2\alpha-3-\frac{1}{\theta}}{1-\frac{n}{2}+\frac{n(m+p-1)}{2}}
    +\frac{\frac{2}{s}-1+\frac{1}{\theta}}{1-\frac{n}{2} + \frac{nq}{s}}
    \qquad \mbox{for $q\ge 2$ and $s \in [1,\frac{n}{n-1}]$},
  \ee
  and
  \be{33.999}
    g(q,s):=\frac{2-\frac{1}{\mu}}{1-\frac{n}{2}+\frac{n(m+p-1)}{2}}
    +\frac{\frac{2(q-1)}{s}-1+\frac{1}{\mu}}{1-\frac{n}{2} + \frac{nq}{s}}
    \qquad \mbox{for $q\ge 2$ and $s \in [1,\frac{n}{n-1}]$}.
  \ee
  Then
  \bas
    g \Big( q_0(p),\frac{n}{n-1} \Big)
    &=& \frac{2-\frac{1}{\mu} + \frac{2(n-1)}{n} \cdot \Big( \frac{n(m+p-1)}{2(n-1)}-1\Big) -1 + \frac{1}{\mu}}
    {1-\frac{n}{2}+\frac{n(m+p-1)}{2}} \\
    &=& \frac{1+(m+p-1)-\frac{2(n-1)}{n}}
    {1-\frac{n}{2}+\frac{n(m+p-1)}{2}} \\
    &=& \frac{\frac{2}{n}-1+(m+p-1)}
    {1-\frac{n}{2}+\frac{n(m+p-1)}{2}} \\
    &=& \frac{2}{n}.
  \eas
  Since
  \bas
    \frac{\partial g}{\partial q} \Big( q, \frac{n}{n-1})
    &=& \frac{\frac{2(n-1)}{n} \cdot \Big[ 1-\frac{n}{2}+(n-1)q \Big]
    - \Big[ \frac{2(n-1)}{n} \cdot (q-1) - 1 + \frac{1}{\mu} \Big] \cdot (n-1)}
    {\Big[ 1-\frac{n}{2} + (n-1)q \Big]^2} \\
    &=& (n-1) \cdot \frac{\frac{2}{n}-1+\frac{2(n-1)q}{n} - \frac{2(n-1)q}{n} + \frac{2(n-1)}{n}+1-\frac{1}{\mu}}
    {\Big[ 1-\frac{n}{2} + (n-1)q \Big]^2} \\
    &=& (n-1) \cdot \frac{2-\frac{1}{\mu}}
    {\Big[ 1-\frac{n}{2} + (n-1)q \Big]^2} \\
    &>& 0 \qquad \mbox{for all } q>2,
  \eas
  this implies
  \be{33.10}
    g \Big(q,\frac{n}{n-1} \Big) < \frac{2}{n}
    \qquad \mbox{for all } q \in (2,q_0(p)).
  \ee
  Moreover, our assumption $\alpha<\frac{2}{n}$ entails that
  \bas
    f \Big( q_0(p), \frac{n}{n-1} \Big)
    &=& \frac{\Big(m+p+2\alpha-3-\frac{1}{\theta}\Big) + \Big( 1-\frac{2}{n}+\frac{1}{\theta}\Big)}
    {1-\frac{n}{2}+\frac{n(m+p-1)}{2}} \\
    &=& \frac{m+p+2\alpha-2-\frac{2}{n}}
    {1-\frac{n}{2}+\frac{n(m+p-1)}{2}} \\
    &<& \frac{m+p-2+\frac{2}{n}}
    {1-\frac{n}{2}+\frac{n(m+p-1)}{2}} \\
    &=& \frac{2}{n} \cdot \frac{\frac{n(m+p-1)}{2}+1-\frac{n}{2}}
    {1-\frac{n}{2}+\frac{n(m+p-1)}{2}} \\
    &=& \frac{2}{n}.
  \eas
  Therefore by a continuity argument using (\ref{33.44}) we can now fix $q>\bar q$ fulfilling
  \be{33.11}
    q<q_0(p)
  \ee
  and
  \be{33.12}
    f \Big(q,\frac{n}{n-1}\Big) < \frac{2}{n}
  \ee
  and such that furthermore (\ref{33.13}) and (\ref{33.14}) hold, where the latter two can be achieved
  on choosing $q$ close enough to $q_0(p)$ according to (\ref{33.8}) and (\ref{33.9}). We observe that by (\ref{33.11})
  and (\ref{33.10}) we also have
  \bas
    g \Big(q,\frac{n}{n-1} \Big)<\frac{2}{n},
  \eas
  so that, again by continuity, we can finally find $s \in [1,\frac{n}{n-1})$ close to $\frac{n}{n-1}$ such that with $q$
  as fixed above we still have
  \bas
    f(q,s)<\frac{2}{n} \qquad \mbox{and} \qquad g(q,s)<\frac{2}{n}.
  \eas
  In view of the definitions (\ref{33.99}) and (\ref{33.999}) of $f$ and $g$, these two inequalities are equivalent
  to (\ref{33.1551}) and (\ref{33.1552}).
\qed
\mysection{Proof of the main results}\label{mainproof}
The following final preparation is a direct consequence of Young's inequality.
\begin{lem}\label{lem_young2}
  Let $\beta>0$ and $\gamma>0$ be such that $\beta+\gamma<1$. Then for all $\eps>0$ there exists $c>0$ such that
  \bas
    a^\beta b^\gamma \le \eps (a+b) + c \qquad \mbox{for all $a \ge 0$ and } b\ge 0.
  \eas
\end{lem}
We proceed to establish the main step towards our boundedness proof.
\begin{lem}\label{lem33}
  Suppose that $\Omega$ is convex, and that
  (\ref{alpha}), (\ref{m}) and (\ref{M}) hold with some $K>0$, $K_0>0, K_1>0, m\in\R, M\in\R$ and some positive
  \bas
    \alpha<\frac{2}{n}.
  \eas
  Then for all $p\in [1,\infty)$ and each $q\in [1,\infty)$ there exists $c>0$ such that
  \be{33.1}
    \|u(t)\|_{L^p(\Omega)} \le c \qquad \mbox{for all } t\in (0,\tm)
  \ee
  and
  \be{33.2}
    \|\nabla v(t)\|_{L^{2q}(\Omega)} \le c \qquad \mbox{for all } t\in (0,\tm)
  \ee
\end{lem}
\proof
  It is evidently sufficient to prove that for any $p_0>1$ and $q_0>2$ we can find some $p>p_0$ and $q>q_0$ such that
  (\ref{33.2}) and
  \be{33.1111}
    \|u(t)\|_{L^{p+m-M}(\Omega)} \le c \qquad \mbox{for all } t\in (0,\tm)
  \ee
  hold with some $c>0$, where $m$ and $M$ are taken from (\ref{m}) and (\ref{M}), respectively.
  To achieve this, given such $p_0$ and $q_0$ let us set $\bar p:=p_0+M-m$ and $\bar q:=q_0$ and then fix
  $p>\bar p, q>\bar q, s \in [1,\frac{n}{n-1}), \theta>1$ and $\mu>1$ as provided by Lemma \ref{lem45}. Then by (\ref{m}),
  \be{33.1555}
    \phi(r) := \int_0^r \int_0^\rho \frac{(\sigma+1)^{m+p-3}}{D(\sigma)} d\sigma d\rho \qquad \mbox{for } r \ge 0.
  \ee
  is finite and positive for all $r\ge 0$ with
  \be{33.1556}
    \phi(r) \le \frac{1}{K_0} \cdot \int_0^r \int_0^\rho (\sigma+1)^{p-2}d\sigma d\rho \le \frac{1}{p(p-1)K_0} \cdot (r+1)^p
    \qquad \mbox{for all } r \ge 0,
  \ee
  and furthermore due to (\ref{M}) we have
  \be{33.1557}
    \phi(r) \ge c_0 (r+1)^{p+m-M} \qquad \mbox{for all } r\ge 0
  \ee
  with some $c_0>0$.
  Since moreover $\phi$ is smooth on $(0,\infty)$ and $u$ is positive in $\Omega \times (0,\tm)$ by the strong
  maximum principle, we may use $\phi'(u)$ as a test function for the first equation in (\ref{0}). Integrating by parts
  we thereby see that
  \bea{33.16}
    \hspace*{-5mm}
    \frac{d}{dt} \io \phi(u)
    &=& \io \phi'(u) \nabla \cdot (D(u)\nabla u) - \io \phi'(u) \nabla \cdot (S(u)\nabla v) \nn\\
    &=& - \io \phi''(u) D(u) |\nabla u|^2 + \io \phi''(u) S(u) \nabla u \cdot \nabla v \nn\\
    &=& - \io (u+1)^{m+p-3}|\nabla u|^2 + \io (u+1)^{m+p-3} \frac{S(u)}{D(u)} \nabla u \cdot \nabla v
  \eea
  for all $t\in (0,\tm)$,
  where thanks to Young's inequality and (\ref{alpha}),
  \be{33.17}
    \io (u+1)^{m+p-3}\frac{S(u)}{D(u)} \nabla u \cdot \nabla v
    \le \frac{1}{2} \io (u+1)^{m+p-3}|\nabla u|^2
    + \frac{K^2}{2} \io (u+1)^{m+p+2\alpha-3}|\nabla v|^2.
  \ee
  We next differentiate the second equation in (\ref{0}) to obtain
  \bas
    (|\nabla v|^2)_t = 2\nabla v \cdot \nabla \Delta v - 2|\nabla v|^2 + 2\nabla u \cdot \nabla v
  \eas
  and hence, recalling the identity $\Delta |\nabla v|^2 = 2\nabla v \cdot \nabla \Delta v + 2|D^2 v|^2$,
  \bas
    (|\nabla v|^2)_t =\Delta |\nabla v|^2 - 2|D^2 v|^2 - 2|\nabla v|^2 + 2\nabla u \cdot \nabla v
  \eas
  for all $x\in \Omega$ and $t\in (0,\tm)$. Testing this against $|\nabla v|^{2q-2}$ yields
  \bea{33.18}
    \frac{1}{q} \frac{d}{dt} \io |\nabla v|^{2q}
    &+& (q-1) \io |\nabla v|^{2q-4} \Big| \nabla |\nabla v|^2 \Big|^2
    +2 \io |\nabla v|^{2q-2} |D^2 v|^2
    + 2\io |\nabla v|^{2q} \nn\\
    &\le& 2 \io |\nabla v|^{2q-2} \nabla u \cdot \nabla v
    \qquad \mbox{for all } t\in (0,\tm),
  \eea
  where we have used the convexity of $\Omega$ which in conjunction with the boundary condition
  $\frac{\partial v}{\partial \nu}=0$ on $\pO$ implies that
  \bas
    \frac{\partial |\nabla v|^2}{\partial \nu} \le 0 \qquad \mbox{on } \pO
  \eas
  (cf.~\cite{dalpasso_garcke_gruen}). On the right of (\ref{33.18}) we integrate by parts and use Young's inequality to find
  \bea{33.19}
    2 \io |\nabla v|^{2q-2} \nabla u \cdot \nabla v
    &=& - 2(q-1) \io u|\nabla v|^{2q-4} \nabla v \cdot \nabla |\nabla v|^2
    - 2 \io u|\nabla v|^{2q-2} \Delta v \nn\\
    &\le& \frac{q-1}{2} \io |\nabla v|^{2q-4} \Big| \nabla |\nabla v|^2 \Big|^2
    + 2(q-1) \io u^2 |\nabla v|^{2q-2} \nn\\
    & & + \frac{2}{n} \io |\nabla v|^{2q-2} |\Delta v|^2
    +\frac{n}{2} \io u^2 |\nabla v|^{2q-2},
  \eea
  where
  \bas
    \frac{2}{n} \io |\nabla v|^{2q-2} |\Delta v|^2 \le 2\io |\nabla v|^{2q-2} |D^2 v|^2
  \eas
  in view of the pointwise inequality $|\Delta v|^2 \le n |D^2 v|^2$.
  We thus infer from (\ref{33.16})-(\ref{33.19}) that there exists $c_1>0$ such that
  \bea{33.20}
    \frac{d}{dt} \bigg\{ \io \phi(u) + \frac{1}{q} \io |\nabla v|^{2q} \bigg\}
    &+& \frac{2}{(m+p-1)^2} \io |\nabla (u+1)^\frac{m+p-1}{2}|^2
    + \frac{2(q-1)}{q^2} \io \Big| \nabla |\nabla v|^q \Big|^2 \nn\\
    &\le& c_1 \io (u+1)^{m+p+2\alpha-3}|\nabla v|^2 + c_1 \io (u+1)^2 |\nabla v|^{2q-2}
  \eea
  for all $t\in (0,\tm)$.
  Here we use the H\"older inequality to estimate the integrals on the right according to
  \be{33.21}
    \io (u+1)^{m+p+2\alpha-3}|\nabla v|^2 \le \Big( \io (u+1)^{(m+p+2\alpha-3)\theta}\Big)^\frac{1}{\theta} \cdot
    \Big( \io |\nabla v|^{2\theta'} \Big)^\frac{1}{\theta'}
  \ee
  and
  \be{33.22}
    \io (u+1)^2 |\nabla v|^{2q-2} \le \Big( \io (u+1)^{2\mu} \Big)^\frac{1}{\mu} \cdot
    \Big( \io |\nabla v|^{2(q-1)\mu'} \Big)^\frac{1}{\mu'}
  \ee
  with $\theta':=\frac{\theta}{\theta-1}$ and $\mu':=\frac{\mu}{\mu-1}$.
  Now since (\ref{33.7}) in conjunction with the positivity of $\alpha$ and the fact that $\theta>1$ implies that
  \bas
    \frac{2(m+p+2\alpha-3)\theta}{m+p-1} > \frac{2}{m+p-1},
  \eas
  and since (\ref{33.5}) asserts that
  \bas
    \frac{2(m+p+2\alpha-3)\theta}{m+p-1} < \frac{2n}{n-2},
  \eas
  we may invoke the Gagliardo-Nirenberg inequality to estimate
  \bea{33.23}
    \Big(\io (u+1)^{(m+p+2\alpha-3)\theta}\Big)^\frac{1}{\theta}
    &=& \|(u+1)^\frac{m+p-1}{2}\|_{L^\frac{2(m+p+2\alpha-3)\theta}{m+p-1}(\Omega)}^\frac{2(m+p+2\alpha-3)}{m+p-1} \nn\\
    & & \hspace*{-20mm}
    \le c_2 \|\nabla (u+1)^\frac{m+p-1}{2}\|_{L^2(\Omega)}^{\frac{2(m+p+2\alpha-3)}{m+p-1} \cdot a}
    \cdot \|(u+1)^\frac{m+p-1}{2}\|_{L^\frac{2}{m+p-1}(\Omega)}^{\frac{2(m+p+2\alpha-3)}{m+p-1} \cdot (1-a)} \nn\\
    & & \hspace*{-20mm}
    + c_2 \|(u+1)^\frac{m+p-1}{2}\|_{L^\frac{2}{m+p-1}(\Omega)}^\frac{2(m+p+2\alpha-3)}{m+p-1}
    \qquad \mbox{for all } t\in (0,\tm),
  \eea
  with some $c_2>0$ and $a\in (0,1)$ determined by
  \bas
    -\frac{n(m+p-1)}{2(m+p+2\alpha-3)\theta}= \Big(1-\frac{n}{2}\Big) \cdot a
    - \frac{n(m+p-1)}{2} \cdot (1-a).
  \eas
  Thus,
  \bas
    a=\frac{\frac{n(m+p-1)}{2} \cdot \big(1-\frac{1}{(m+p+2\alpha-3)\theta}\big)}
    {1-\frac{n}{2}+\frac{n(m+p-1)}{2}}
  \eas
  and hence
  \bas
    \frac{2(m+p+2\alpha-3)}{m+p-1} \cdot a =
    n \cdot \frac{m+p+2\alpha-3-\frac{1}{\theta}}{1-\frac{n}{2}+\frac{n(m+p-1)}{2}},
  \eas
  so that (\ref{33.23}) yields
  \be{33.24}
    \Big(\io (u+1)^{(m+p+2\alpha-3)\theta} \Big)^\frac{1}{\theta}
    \le c_3 \Big( \io |\nabla (u+1)^\frac{m+p-1}{2}|^2 \Big)^
    {\frac{n}{2} \cdot \frac{m+p+2\alpha-3-\frac{1}{\theta}}{1-\frac{n}{2}+\frac{n(m+p-1)}{2}}}  + c_3
  \ee
  for all $t\in (0,\tm)$ with some $c_3>0$, because (\ref{mass}) states boundedness of $(u+1)^\frac{m+p-1}{2}$ in
  $L^\infty((0,\tm);L^\frac{2}{m+p-1}(\Omega))$. \abs
  Similarly, using that $\mu>1$ implies
  \bas
    \frac{4\mu}{m+p-1}>\frac{2}{m+p-1},
  \eas
  and that (\ref{33.6}) entails
  \bas
    \frac{4\mu}{m+p-1} < \frac{2n}{n-2},
  \eas
  we interpolate
  \bas
    \Big( \io (u+1)^{2\mu} \Big)^\frac{1}{\mu}
    &=& \|(u+1)^\frac{m+p-1}{2}\|_{L^\frac{4\mu}{m+p-1}(\Omega)}^\frac{4}{m+p-1} \\
    &\le& c_4 \|\nabla (u+1)^\frac{m+p-1}{2} \|_{L^2(\Omega)}^{\frac{4}{m+p-1} \cdot b}
    \cdot \|(u+1)^\frac{m+p-1}{2} \|_{L^\frac{2}{m+p-1}(\Omega)}^{\frac{4}{m+p-1} \cdot (1-b)} \\
    & & + c_4 \|(u+1)^\frac{m+p-1}{2}\|_{L^\frac{2}{m+p-1}(\Omega)}^\frac{4}{m+p-1}
    \qquad \mbox{for all } t\in (0,\tm)
  \eas
  with some $c_4>0$ and
  \bas
    b=\frac{\frac{n(m+p-1)}{2} \big(1-\frac{1}{2\mu}\big)}{1-\frac{n}{2}+\frac{n(m+p-1)}{2}} \, \in (0,1).
  \eas
  Again in view of (\ref{mass}), we therefore obtain $c_5>0$ such that
  \be{33.25}
    \Big(\io (u+1)^{2\mu} \Big)^\frac{1}{\mu}
    \le c_5 \Big( \io |\nabla (u+1)^\frac{m+p-1}{2}|^2 \Big)^
    {\frac{n}{2} \cdot \frac{2-\frac{1}{\mu}}{1-\frac{n}{2}+\frac{n(m+p-1)}{2}}}    + c_5
    \qquad \mbox{for all } t\in (0,\tm).
  \ee
  As to the integrals in (\ref{33.21}) and (\ref{33.22}) involving $\nabla v$, we proceed in quite the same manner,
  relying on (\ref{gradv}) rather than on (\ref{mass}). First, we note that
  \be{33.26}
    \frac{2\theta'}{q} > \frac{s}{q},
  \ee
  because $\theta'>1$ and $s<\frac{n}{n-1} \le 2$ whenever $n\ge 2$. Moreover, we know that
  \be{33.27}
    \frac{2\theta'}{q} < \frac{2n}{n-2},
  \ee
  for (\ref{33.13}) says that
  \bas
    \frac{1}{\theta'} = 1-\frac{1}{\theta}>\frac{n-2}{n} \cdot \frac{1}{q}.
  \eas
  Now (\ref{33.26}) and (\ref{33.27}) allow for an application of the Gagliardo-Nirenberg inequality which ensures the
  existence of $c_6>0$ fulfilling
  \bas
    \Big( \io |\nabla v|^{2\theta'} \Big)^\frac{1}{\theta'}
    &=& \Big\| |\nabla v|^q \Big\|_{L^\frac{2\theta'}{q}(\Omega)}^\frac{2}{q} \\
    &\le& c_6 \Big\| \nabla |\nabla v|^q \Big\|_{L^2(\Omega)}^{\frac{2}{q} \cdot c} \cdot
    \Big\| |\nabla v|^q \Big\|_{L^\frac{s}{q}(\Omega)}^{\frac{2}{q}(1-c)}
    + c_6 \Big\| |\nabla v|^q \Big\|_{L^\frac{s}{q}(\Omega)}^\frac{2}{q}
  \eas
  with
  \bas
    c= \frac{nq(\frac{1}{s}-\frac{1}{2\theta'})}{1-\frac{n}{2}+\frac{nq}{s}} \, \in (0,1).
  \eas
  By means of (\ref{gradv}), we thus find $c_7>0$ such that
  \be{33.28}
    \Big(\io |\nabla v|^{2\theta'} \Big)^\frac{1}{\theta'}
    \le c_7 \bigg( \io \Big|\nabla |\nabla v|^q \Big|^2 \bigg)^
    {\frac{n}{2} \cdot \frac{\frac{2}{s}-\frac{1}{\theta'}}{1-\frac{n}{2}+\frac{nq}{s}}}    + c_7
    \qquad \mbox{for all } t\in (0,\tm).
  \ee
  As to the corresponding term in (\ref{33.22}), we similarly observe that
  \be{33.29}
    \frac{2(q-1)\mu'}{q} > \frac{s}{q},
  \ee
  which immediately follows from the inequalities $\mu'>1$ and $q>\bar q\ge 2$ and our assumption $n\ge 2$.
  We furthermore have
  \be{33.30}
    \frac{2(q-1)\mu'}{q} < \frac{2n}{n-2},
  \ee
  because (\ref{33.14}) asserts that
  \bas
    \frac{1}{\mu'} = 1-\frac{1}{\mu} > 1-\frac{2}{n} -\frac{n-2}{n} \cdot \frac{1}{q} = \frac{n-2}{n} \cdot \frac{q-1}{q}.
  \eas
  Thanks to (\ref{33.29}), (\ref{33.30}) and the Gagliardo-Nirenberg inequality, we can find $c_8>0$ satisfying
  \bas
    \Big( \io |\nabla v|^{2(q-1)\mu'} \Big)^\frac{1}{\mu'}
    &=& \Big\| |\nabla v|^q \Big\|_{L^\frac{2(q-1)\mu'}{q}(\Omega)}^\frac{2(q-1)}{q} \\
    &\le& c_8 \Big\| \nabla |\nabla v|^q \Big\|_{L^2(\Omega)}^{\frac{2(q-1)}{q} \cdot d} \cdot
    \Big\| |\nabla v|^q \Big\|_{L^\frac{s}{q}(\Omega)}^{\frac{2(q-1)}{q} \cdot (1-d)}
    + c_8 \Big\| æ\nabla v|^q \big\|_{L^\frac{s}{q}(\Omega)}^\frac{2(q-1)}{q}
  \eas
  with
  \bas
    d=\frac{nq \cdot (\frac{1}{s} - \frac{1}{2(q-1)\mu'})}{1-\frac{n}{2}+\frac{nq}{s}} \, \in (0,1).
  \eas
  Consequently, once again recalling (\ref{gradv}) we have
  \be{33.31}
    \Big(\io |\nabla v|^{2(q-1)\mu'}\Big)^\frac{1}{\mu'}
    \le c_9 \bigg( \io \Big| \nabla |\nabla v|^q \Big|^2 \bigg)^
    {\frac{n}{2} \cdot \frac{\frac{2(q-1)}{s} - \frac{1}{\mu'}}{1-\frac{n}{2}+\frac{nq}{s}}}    + c_9
    \qquad \mbox{for all } t\in (0,\tm)
  \ee
  for some positive constant $c_9$.\abs
  Now collecting (\ref{33.24}), (\ref{33.25}), (\ref{33.28}) and (\ref{33.31}), from (\ref{33.21}) and (\ref{33.22}) we obtain
  \bea{33.32}
    & & \hspace*{-30mm}
    c_1 \io (u+1)^{m+p+2\alpha-3}|\nabla v|^2 + c_1 \io (u+1)^2 |\nabla v|^{2q-2} \nn\\
    &\le& c_{10} \bigg( \io |\nabla (u+1)^\frac{m+p-1}{2}|^2 \bigg)^{\beta_1} \cdot
    \bigg( \io \Big| \nabla |\nabla v|^q \Big|^2 \bigg)^{\gamma_1} \nn\\
    & & + c_{10} \bigg( \io |\nabla (u+1)^\frac{m+p-1}{2}|^2 \bigg)^{\beta_2} \cdot
    \bigg( \io \Big| \nabla |\nabla v|^q \Big|^2 \bigg)^{\gamma_2} \nn\\[1mm]
    & & + c_{10}
    \qquad \mbox{for all } t\in (0,\tm)
  \eea
  with some $c_{10}>0$ and positive numbers $\beta_1, \beta_2,\gamma_1$ and $\gamma_2$ satisfying
  \bas
    \beta_1+\gamma_1 &=& \frac{n}{2} \cdot \frac{m+p+2\alpha-3-\frac{1}{\theta}}{1-\frac{n}{2}+\frac{n(m+p-1)}{2}}
    +\frac{n}{2} \cdot \frac{\frac{2}{s}-1+\frac{1}{\theta}}{1-\frac{n}{2}+\frac{nq}{s}} \\
    &<& 1
  \eas
  according to (\ref{33.1551}), and
  \bas
    \beta_2+\gamma_2 &=& \frac{n}{2} \cdot \frac{2-\frac{1}{\mu}}{1-\frac{n}{2}+\frac{n(m+p-1)}{2}}
    +\frac{n}{2} \cdot \frac{\frac{2(q-1)}{s}-1+\frac{1}{\mu}}{1-\frac{n}{2}+\frac{nq}{s}} \\
    &<& 1
  \eas
  by (\ref{33.1552}).
  Therefore Lemma \ref{lem_young2} states that for some $c_{11}>0$ we have
  \bea{33.33}
    & & \hspace*{-20mm}
    c_1 \io (u+1)^{m+p+2\alpha-3}|\nabla v|^2 + c_1 \io (u+1)^2 |\nabla v|^{2q-2} \nn\\
    &\le& \frac{1}{(m+p-1)^2} \io |\nabla (u+1)^\frac{m+p-1}{2}|^2
    + \frac{q-1}{q^2} \io \Big| \nabla |\nabla v|^q \Big|^2  + c_{11}
  \eea
  for all $t\in (0,\tm)$. Here we once more employ the Gagliardo-Nirenberg inequality to estimate
  \bea{33.34}
    \hspace*{-5mm}
    \io (u+1)^p &=& \|(u+1)^\frac{m+p-1}{2} \|_{L^\frac{2p}{m+p-1}(\Omega)}^\frac{2p}{m+p-1} \nn\\
    &\le& c_{12} \|\nabla (u+1)^\frac{m+p-1}{2}\|_{L^2(\Omega)}^{\frac{2p}{m+p-1} \cdot \kappa_1}
    \cdot \|(u+1)^\frac{m+p-1}{2}\|_{L^\frac{2}{m+p-1}(\Omega)}^{\frac{2p}{m+p-1} \cdot (1-\kappa_1)} \nn\\
    & & + c_{12} \|(u+1)^\frac{m+p-1}{2}\|_{L^\frac{2}{m+p-1}(\Omega)}^{\frac{2p}{m+p-1}}
  \eea
  and
  \bea{33.35}
    \io |\nabla v|^{2q} &=& \Big\| |\nabla v|^q \Big\|_{L^2(\Omega)}^2 \nn\\
    &\le& c_{12} \Big\| \nabla |\nabla v|^q \Big\|_{L^2(\Omega)}^{2\kappa_2}
    \cdot \Big\| |\nabla v|^q \Big\|_{L^\frac{s}{q}(\Omega)}^{2(1-\kappa_2)}
    + c_{12} \Big\| |\nabla v|^q \Big\|_{L^\frac{s}{q}(\Omega)}^2
  \eea
  with some $c_{12}>0$ and
  \bas
    \kappa_1=\frac{\frac{n(m+p-1)}{2}(1-\frac{1}{p})}{1-\frac{n}{2}+\frac{n(m+p-1)}{2}}
    \qquad \mbox{and} \qquad
    \kappa_2=\frac{\frac{nq}{s}-\frac{n}{2}}{1-\frac{n}{2}+\frac{nq}{s}},
  \eas
  where we note that $\frac{2p}{m+p-1}<\frac{2n}{n-2}$ by (\ref{33.7}) and $\frac{s}{q}<2$ since $q>\bar q > 2$ and
  $s<\frac{n}{n-1} \le 2$.\\
  As a consequence of (\ref{33.34}), (\ref{33.35}), (\ref{mass}) and (\ref{gradv}), (\ref{33.20}) can be turned into
  the inequality
  \bas
    \frac{d}{dt} \bigg( \io \phi(u) + \frac{1}{q} \io |\nabla v|^{2q} \bigg)
    +c_{13} \Big( \io (u+1)^p \Big)^\frac{m+p-1}{p\kappa_1} + c_{13} \Big( \io |\nabla v|^{2q} \Big)^\frac{1}{\kappa_2}
    \le c_{14}
  \eas
  for all $t\in (0,\tm)$ and positive constants $c_{13}$ and $c_{14}$. In view of (\ref{33.1556}), we infer that the
  function
  \bas
    y(t):= \io \phi(u(t)) + \frac{1}{q} \io |\nabla v(t)|^{2q}, \qquad t\in [0,\tm),
  \eas
  satisfies
  \bas
    y'(t) + c_{15} y^\kappa(t) \le c_{16}
    \qquad \mbox{for all } t\in (0,\tm)
  \eas
  with certain positive constants $\kappa, c_{15}$ and $c_{16}$. Upon an ODE comparison argument this entails that
  \bas
    y(t) \le c_{17} := \max \bigg\{ y_0, \Big( \frac{c_{16}}{c_{15}} \Big)^\frac{1}{\kappa} \bigg\}
    \qquad \mbox{for all } t\in (0,\tm).
  \eas
  Thus, in view of (\ref{33.1557}) we arrive at the inequalities
  \bas
    \io (u+1)^{p+m-M}(t) \le \frac{c_{17}}{c_0}
    \qquad \mbox{and} \qquad
    \io |\nabla v(t)|^{2q} \le qc_{17}
    \qquad \mbox{for all } t\in (0,\tm)
  \eas
  and conclude.
\qed
Now we can immediately pass to our main result.\abs
\proofc (of Theorem \ref{theo_main}) \quad
  The proof is an evident consequence of Lemma \ref{lem33}, Lemma \ref{lem_moser} below and the
  extendibility criterion provided by Lemma \ref{lem_loc_exist}.
\qed
\mysection{Appendix: A general boundedness result for quasilinear non-uniformly parabolic equations}
In this concluding section, which might be of interest of its own, we derive uniform bounds for nonnegative subsolutions of some
quasilinear problems which need not necessarily be uniformly parabolic.
More precisely, we consider functions $u$ fulfilling
\be{m1}
    \left\{ \begin{array}{l}
    u_t \le \nabla \cdot (D(x,t,u)\nabla u) + \nabla \cdot f(x,t) + g(x,t), \qquad x\in \Omega, \ t\in (0,T), \\[1mm]
    \dN u(x,t) \le 0, \qquad x\in\pO, \ t \in (0,T),
    \end{array} \right.
\ee
in the classical sense, where we allow the diffusion to be degenerate in the sense that we require that
\be{D1}
    D \in C^1(\bar\Omega \times [0,T) \times [0,\infty)) \qquad \mbox{and} \qquad
    D \ge 0,
\ee
and that there exist $m\in\R$, $s_0 \ge 0$ and $\delta>0$ such that
\be{D}
    D(x,t,s) \ge \delta s^{m-1} \qquad \mbox{for all $x\in\Omega$, $t \in (0,T)$ and } s \ge s_0.
\ee
%
%
%
Our goal is to derive estimates in $L^\infty(\Omega \times (0,T))$ under the assumptions
that
\be{reg_fg}
    f\in C^0((0,T);C^0(\bar\Omega) \cap C^1(\Omega)) \qquad \mbox{and} \qquad
    g \in C^0(\Omega \times (0,T))
\ee
with
\be{bdry_f}
    f\cdot \nu \le 0 \qquad \mbox{on } \pO \times (0,T),
\ee
that
\be{fg}
    f\in L^\infty((0,T);L^{q_1}(\Omega)) \qquad \mbox{and} \qquad g \in L^\infty((0,T);L^{q_2}(\Omega)),
\ee
and that
\be{u_infty}
    u \in L^\infty((0,T);L^{p_0}(\Omega))
\ee
be valid with suitably large $q_1, q_2$ and $p_0$.\\
The derivation of the following statement follows a well-established iterative technique (see \cite{alikakos}
for an application in a similar framework). Since we could not find a precise reference covering our situation,
and since some major modifications to the original procedure are necessary, we inculde a full proof here for the sake of
completeness.
\begin{lem}\label{lem_moser}
  Suppose that $T \in (0,\infty]$, that $\Omega \subset \R^n$ is a bounded domain, and that $D, f$ and $g$ comply with
  (\ref{D1}), (\ref{reg_fg}) and (\ref{bdry_f}).
  Moreover, assume that (\ref{D}) and (\ref{fg}) hold for some $\delta>0$, $m\in\R$ and $s_0 \ge 0$, and some
  $q_1>n+2$ and $q_2>\frac{n+2}{2}$.
  Then if $u \in C^0(\bar\Omega \times [0,T)) \cap C^{2,1}(\bar\Omega \times (0,T))$ is a nonnegative function
  satisfying (\ref{m1}), and if (\ref{u_infty})
  is valid for some $p_0\ge 1$ fulfilling
  \be{p1}
    p_0 > 1-m\cdot \frac{(n+1)q_1-(n+2)}{q_1-(n+2)}
  \ee
  and
  \be{p2}
    p_0 > 1-\frac{m}{1-\frac{n}{n+2}\frac{q_2}{q_2-1}}
  \ee
  as well as
  \be{pp}
    p_0 > \frac{n(1-m)}{2},
  \ee
  then there exists $C>0$, only depending on $m,\delta,\Omega,
  \|f\|_{L^\infty((0,T);L^{q_1}(\Omega))}, \|g\|_{L^\infty((0,T);L^{q_2}(\Omega))}$,
  $\|u\|_{L^\infty((0,T);L^{p_0}(\Omega))}$ and $\|u(0)\|_{L^\infty(\Omega)}$, such that
  \bas
    \|u(t)\|_{L^\infty(\Omega)} \le C \qquad \mbox{for all } t\in (0,T).
  \eas
\end{lem}
\proof
  We evidently may assume that $m \le 0$,
  and then fix $r \in (2,\frac{2(n+2)}{n})$ close enough to $\frac{2(n+2)}{n}$ such that
  writing $\theta(\rho):=\frac{\rho}{2} \cdot \frac{m+p_0-1}{-m+p_0-1}$ and $\mu(\rho):=\frac{\rho}{2} \cdot
  \frac{m+p_0-1}{p_0-1}$ we have $\theta(r) \ge \frac{q_1}{q_1-2}$ and $\mu(r) \ge \frac{q_2}{q_2-1}$. Indeed, this is possible
  since our assumption (\ref{p1}) on $p_0$ ensures that
  \bas
    \theta \Big( \frac{2(n+2)}{n} \Big) &=& \frac{n+2}{n} \cdot \Big( 1+ \frac{2m}{-m+p_0-1} \Big) \\
    &>& \frac{n+2}{n} \cdot \Bigg( 1 + \frac{2m}{-m+[1-m\cdot \frac{(n+1)q_1-(n+2)}{q_1-(n+2)}]-1} \Bigg) \\
    &=& \frac{n+2}{n}\cdot \frac{nq_1}{(n+2)(q_1-2)} \\
    &=& \frac{q_1}{q_1-2}
  \eas
  due to the fact that $q_1>n+2$, and since (\ref{p2}) entails
  \bas
    \mu \Big(\frac{2(n+2)}{n} \Big) &=& \frac{n+2}{n} \cdot \Big(1+\frac{m}{p_0-1}\Big) \\
    &>& \frac{n+2}{n} \cdot \Bigg( 1+\frac{m}{\Big[ 1-\frac{m}{1-\frac{n}{n+2} \cdot \frac{q_2}{q_2-1}} \Big]-1} \Bigg) \\
    &=& \frac{q_2}{q_2-1}.
  \eas
  We can now pick $s\in (0,2)$ sufficiently close to $2$ fulfilling
  \be{m3}
    r<\frac{2(n+s)}{n}
  \ee
  and such that
  \be{m33}
    \frac{\frac{nr}{s}-n}{\frac{2q_1}{q_1-2} \cdot (1-\frac{n}{2}+\frac{n}{s})} < 1,
  \ee
  where the latter can be achieved due to the fact that as $s\to 2$, the expression on the left tends to
  \bas
    \frac{\frac{nr}{2}-n}{\frac{2q_1}{q_1-2}} < \frac{\frac{n}{2} \cdot \frac{2(n+2)}{n} - n}{\frac{2q_1}{q_1-2}}
    =1-\frac{2}{q_1}<1.
  \eas
  We now recursively define
  \be{m4}
    p_k:=\frac{2}{s} \cdot p_{k-1} + 1 -m, \qquad k\ge 1,
  \ee
  and note that $(p_k)_{k\in\N}$ increases and
  \be{m5}
    c_1 \cdot \Big(\frac{2}{s} \Big)^k \le p_k \le c_2 \cdot \Big(\frac{2}{s}\Big)^k \qquad \mbox{for all } k\in\N
  \ee
  holds with positive $c_1$ and $c_2$ which, as all constants $c_3,c_4,...$ appearing below, are independent of $k$.
  Writing
  \be{m6}
    \theta_k:=\frac{r}{2} \cdot \frac{m+p_k-1}{-m+p_k-1}, \qquad k\in\N,
  \ee
  since $m\le 0$ we see that also $(\theta_k)_{k\in\N}$ is increasing with $\theta_k \ge \theta_0=\theta(r)\ge \frac{q_1}{q_1-2}$,
  and hence $\theta_k':=\frac{\theta_k}{\theta_k-1}$ satisfies
  \be{m7}
    1<\theta_k' \le \frac{q_1}{2} \qquad \mbox{for all } k\in\N.
  \ee
  Similarly,
  \be{m66}
    \mu_k:=\frac{r}{2} \cdot \frac{m+p_k-1}{p_k-1}, \qquad k\in\N,
  \ee
  defines an increasing sequence of numbers such that $\mu_k \ge \mu_0=\mu(r)\ge \frac{q_2}{q_2-1}$,
  and such that for $\mu_k':=\frac{\mu_k}{\mu_k-1}$ we have
  \be{m777}
    1<\mu_k' \le q_2 \qquad \mbox{for all } k\in\N.
  \ee
  Our goal is to derive a recursive inequality for
  \be{m77}
    M_k:=\sup_{t\in (0,T)} \io \hu^{p_k}(x,t) dx, \qquad k\in\N,
  \ee
  where $\hu(x,t):=\max\{u(x,t),s_0\}$ for $x\in \bar\Omega$ and $t\in [0,T)$.
  To this end, we note that by a standard approximation procedure we may use $p_k \hu^{p_k-1}$ as a test function
  in (\ref{m1}) to obtain for $k\ge 1$
  \bas
    \frac{d}{dt} \io \hu^{p_k} &+& p_k(p_k-1) \io D(x,t,u) \hu^{p_k-2} |\nabla \hu|^2 \\
    &\le& -p_k(p_k-1) \io \hu^{p_k-2} f\cdot \nabla \hu
    + p_k \io \hu^{p_k-1} g
  \eas
  for all $t\in (0,T)$, where we have used our assumptions that $f\cdot \nu \le 0 $ and $\dN u \le 0$ on $\pO$.
  We now employ Young's inequality in estimating
  \bas
    -p_k(p_k-1) \io \hu^{p_k-2} f\cdot \nabla \hu
    &\le& \frac{p_k(p_k-1)\delta}{2} \io \hu^{m+p_k-3} |\nabla \hu|^2 \\
    & & + \frac{p_k(p_k-1)}{2\delta} \io \hu^{-m+p_k-1} |f|^2,
  \eas
  recall (\ref{D}) and observe that $D(x,t,u)=D(x,t,\hu)$ wherever $u\ge s_0$, to find $c_3>0$ and $c_4>0$ such that
  \bea{m8}
    \frac{d}{dt} \io \hu^{p_k} + c_3 \io \Big|\nabla \hu^\frac{m+p_k-1}{2} \Big|^2
    &\le& c_4 p_k^2 \io \hu^{-m+p_k-1} |f|^2 \nn\\
    & & + p_k \io \hu^{p_k-1} g
  \eea
  for all $t\in (0,T)$. Here, by the H\"older inequality, (\ref{fg}) and (\ref{m7}), there exists $c_5>0$ such that
  \bas
    \io \hu^{-m+p_k-1} |f|^2 &\le& \Big( \io \hu^{(-m+p_k-1)\theta_k} \Big)^\frac{1}{\theta_k}
    \cdot \big( \io |f|^{q_1} \Big)^\frac{2}{q_1} \cdot |\Omega|^\frac{q_1-2\theta_k'}{q_1 \theta_k'} \nn\\
    &\le& c_5 \Big\|\hu^\frac{m+p_k-1}{2} \Big\|_
    {L^\frac{2(-m+p_k-1)\theta_k}{m+p_k-1}(\Omega)}^\frac{2(-m+p_k-1)}{m+p_k-1} \nn\\
    &=& c_5 \Big\|\hu^\frac{m+p_k-1}{2}\Big\|_{L^r(\Omega)}^\frac{r}{\theta_k}
    \qquad \mbox{for all } t\in (0,T)
  \eas
  due to (\ref{m6}).
  Similarly, thanks to (\ref{m777}) there exists $c_6>0$ such that
  \bas
    \io \hu^{p_k-1}g
    &\le& \Big( \io \hu^{(p_k-1)\mu_k} \Big)^\frac{1}{\mu_k} \cdot
    \Big( \io |g|^{q_2} \Big)^\frac{1}{q_2} \cdot |\Omega|^\frac{q_2-\mu_k'}{q_2 \mu_k'} \\
    &\le& c_6 \Big\|\hu^\frac{m+p_k-1}{2}\Big\|_{L^r(\Omega)}^\frac{r}{\mu_k}
    \qquad \mbox{for all } t\in (0,T).
  \eas
  Since evidently $\mu_k\ge \theta_k$ and $p_k \ge 1$ for $k\ge 1$, from (\ref{m8}) we thus see that
  \be{m88}
    \frac{d}{dt} \io \hu^{p_k} + c_3 \io \Big|\nabla \hu^\frac{m+p_k-1}{2} \Big|^2
    \le c_7 p_k^2 \Big\|\hu^\frac{m+p_k-1}{2}\Big\|_{L^r(\Omega)}^\frac{r}{\theta_k}
    \qquad \mbox{for all } t\in (0,T)
  \ee
  is valid with some $c_7>0$.\\
  Now invoking the Gagliardo-Nirenberg inequality (\cite{friedman_book}) we find $c_8>0$, by (\ref{m7}) yet independent of $k$,
  such that
  \bas
    \Big\|\hu^\frac{m+p_k-1}{2}\Big\|_{L^r(\Omega)}^\frac{r}{\theta_k}
    &\le& c_8 \Big\|\nabla \hu^\frac{m+p_k-1}{2} \Big\|_{L^2(\Omega)}^\frac{ra}{\theta_k}
    \cdot \Big\| \hu^\frac{m+p_k-1}{2} \Big\|_{L^s(\Omega)}^\frac{r(1-a)}{\theta_k} \\
    & & +c_8 \Big\| \hu^\frac{m+p_k-1}{2} \Big\|_{L^s(\Omega)}^\frac{r}{\theta_k},
  \eas
  whence observing that $\frac{(m+p_k-1)s}{2}=p_{k-1}$ by (\ref{m4}), from (\ref{m77}) we obtain
  \bas
    \Big\| \hu^\frac{m+p_k-1}{2} \Big\|_{L^r(\Omega)}^\frac{r}{\theta_k}
    \le c_8 M_{k-1}^\frac{r(1-a)}{\theta_k s} \cdot
    \bigg( \io \Big| \nabla \hu^\frac{m+p_k-1}{2} \Big|^2 \bigg)^\frac{ra}{2\theta_k}
    +c_8 M_{k-1}^\frac{r}{\theta_k s}
  \eas
  for all $t\in (0,T)$, with
  \be{m99}
    a=\frac{\frac{n}{s}-\frac{n}{r}}{1-\frac{n}{2}+\frac{n}{s}} \, \in (0,1).
  \ee
  Upon an application of Young's inequality, (\ref{m88}) thus yields
  \bea{m888}
    \frac{d}{dt} \io \hu^{p_k} + \frac{c_3}{2} \io \Big|\nabla \hu^\frac{m+p_k-1}{2} \Big|^2
    \le c_9 \Big( p_k^2 M_{k-1}^\frac{r(1-a)}{\theta_k s} \Big)^\frac{2\theta_k}{2\theta_k-ra}
    + c_9 p_k^2 M_{k-1}^\frac{r}{\theta_k s}
  \eea
  for all $t\in (0,T)$ and some $c_9>0$, where we made use of the fact that (\ref{m33}) entails that
  \bas
    \frac{ra}{2\theta_k} \le \frac{ra}{2\theta_0} \le \frac{ra}{\frac{2q_1}{q_1-2}}
    = \frac{\frac{nr}{s}-n}{\frac{2q_1}{q_1-2} \cdot (1-\frac{n}{2}+\frac{n}{s})} < 1
    \qquad \mbox{for all } k\in\N.
  \eas
  Next, since $p_k>\frac{n(1-m)}{2}$ for all $k\ge 1$ by (\ref{pp}), we can pick $\lambda \in (2,\frac{2n}{(n-2)_+})$ such that
  $\frac{2p_k}{m+p_k-1} \le \lambda$ for all $k\ge 1$. Thus, by the H\"older inequality,
  \bas
    \io \hu^{p_k} &=& \Big\| \hu^\frac{m+p_k-1}{2}\Big\|_{L^\frac{2p_k}{m+p_k-1}(\Omega)}^\frac{2p_k}{m+p_k-1} \\
    &\le& |\Omega|^{1-\frac{2p_k}{\lambda (m+p_k-1)}} \cdot \Big\| \hu^\frac{m+p_k-1}{2}\Big\|_
        {L^\lambda(\Omega)}^\frac{2p_k}{m+p_k-1}\\
    &\le& c_{10} \Big\| \hu^\frac{m+p_k-1}{2} \Big\|_{L^\lambda(\Omega)}^\frac{2p_k}{m+p_k-1}
    \qquad \mbox{for all } t\in (0,T)
  \eas
  with some $c_{10}>0$, and therefore applying the Poincar\'e inequality in the form
  \bas
    \|\varphi\|_{L^\lambda(\Omega)}^2 \le c_{11} \Big( \|\nabla \varphi\|_{L^2(\Omega)}^2 + \|\varphi\|_{L^s(\Omega)}^2 \Big)
    \qquad \mbox{for all } \varphi\in W^{1,2}(\Omega),
  \eas
  we infer that
  \bas
    \io \hu^{p_k} \le c_{12} \cdot \bigg( \Big\|\nabla \hu^\frac{m+p_k-1}{2} \Big\|_{L^2(\Omega)}^2
    + \Big\| \hu^\frac{m+p_k-1}{2} \Big\|_{L^s(\Omega)}^2 \bigg)^\frac{p_k}{m+p_k-1}
  \eas
  holds for all $t\in (0,T)$ with a certain $c_{12}>0$.
  In consequence, writing $c_{13}:=\inf_{k\ge 1} c_{12}^{-\frac{m+p_k-1}{p_k}}>0$, we have
  \bas
    \io \Big| \nabla \hu^\frac{m+p_k-1}{2} \Big|^2 \ge c_{13} \Big( \io \hu^{p_k} \Big)^\frac{m+p_k-1}{p_k}
    -M_{k-1}^\frac{2}{s}
    \qquad \mbox{for all } t\in (0,T).
  \eas
  Combined with (\ref{m888}), this gives the inequality
  \bea{m11}
    \frac{d}{dt} \io \hu^{p_k} &\le& -\frac{c_3}{2} \cdot c_{13} \cdot \Big( \io \hu^{p_k} \Big)^\frac{m+p_k-1}{p_k} \nn\\
    & & + c_9 p_k^\frac{4\theta_k}{2\theta_k-ra} \cdot M_{k-1}^\frac{2r(1-a)}{s(2\theta_k-ra)}
    + c_9 p_k^2 M_{k-1}^\frac{r}{\theta_k s}
    + \frac{c_3}{2} M_{k-1}^\frac{2}{s}
  \eea
  for all $t\in (0,T)$ and $k\ge 1$. To simplify this, we observe that
  \bas
    \frac{2r(1-a)}{s(2\theta_k-ra)} \ge \max \Big\{ \frac{r}{\theta_k s}, \frac{2}{s} \Big\} \qquad \mbox{for all }k\ge 1,
  \eas
  because (\ref{m6}) guarantees that $\theta_k \le \frac{r}{2}$. Since furthermore clearly
  \bas
    2<\frac{4\theta_k}{2\theta_k-ra} \le \frac{4\theta_0}{2\theta_0-ra} \qquad \mbox{for all }  k\ge 1,
  \eas
  from (\ref{m11}) and (\ref{m5}) we obtain
  \bas
    \frac{d}{dt} \io \hu^{p_k} &\le& -c_{14} \Big( \io \hu^{p_k} \Big)^\frac{m+p_k-1}{p_k}
    + c_{15} \cdot \tilde b^k \cdot M_{k-1}^\frac{2r(1-a)}{s(2\theta_k-ra)}
  \eas
  for all $t\in (0,T)$ and $k\ge 1$, suitable $c_{14}>0$ and $c_{15}>0$
  and $\tilde b:=(\frac{2}{s})^\frac{4\theta_0}{2\theta_0-ra}>1$.\\
  An integration of this ODI provides $c_{16}>0$ such that
  \bea{m12}
    M_k &\le& \max \Bigg\{ \io \hu_0^{p_k}, \,
    \Big[ \frac{c_{15}}{c_{14}} \cdot \tilde b^k \cdot
    M_{k-1}^\frac{2r(1-a)}{s(2\theta_k-ra)} \Big]^\frac{p_k}{m+p_k-1} \Bigg\} \nn\\
    &\le& \max \bigg\{ \io \hu_0^{p_k}, \, c_{16} b^k M_{k-1}^{\kappa_k} \bigg\}
    \qquad \mbox{for all } k\ge 1,
  \eea
  where $\hu_0(x):=\hu(x,0)$ for $x\in\Omega$, $\kappa_k:=\frac{2r(1-a)}{s(2\theta_k-ra)} \cdot \frac{p_k}{m+p_k-1}$ and
  $b:=\tilde b^\frac{p_0}{m+p_0-1}$, and where we have used that $\frac{p_k}{m+p_k-1} \le \frac{p_0}{m+p_0-1}$ for all
  $k\ge 1$.
  Writing
  \bas
    \kappa_k=\frac{2}{s} \cdot \bigg( 1+\frac{1-\frac{2\theta_k}{r}}{\frac{2\theta_k}{r}-a} \bigg) \cdot
    \bigg( 1+\frac{1-m}{m+p_k-1} \bigg),
  \eas
  we easily infer from (\ref{m4}), (\ref{m6}) and (\ref{m5}) that
  \be{m13}
    \kappa_k = \frac{2}{s} \cdot (1+\eps_k), \qquad k\ge 1,
  \ee
  holds with some $\eps_k \ge 0$ satisfying
  \be{m14}
    \eps_k \le \frac{c_{17}}{p_k} \le c_{18} \cdot \Big( \frac{s}{2} \Big)^k
  \ee
  for all $k\ge 1$ and appropriately large $c_{17}>0$ and $c_{18}>0$. \\
  Therefore, in the case when $c_{16} b^k M_{k-1}^{\kappa_k} < \io \hat{u}_0^{p_k}$ holds for infinitely
  many $k \ge 1$, we obtain
  \bas
    \sup_{t\in (0,T)} \Big( \io \hu^{p_{k-1}} \Big)^\frac{1}{p_{k-1}}
    \le \bigg( \frac{1}{c_{16} b^k} \io \hu_0^{p_k} \bigg)^\frac{1}{\kappa_k p_{k-1}}
  \eas
  for all such $k$, and hence conclude that
  \bas
    \|\hu(t)\|_{L^\infty(\Omega)} \le \|\hu_0\|_{L^\infty(\Omega)} \qquad \mbox{for all } t\in (0,T),
  \eas
  because $\frac{p_k}{\kappa_k p_{k-1}} \to 1$ as $k\to\infty$ according to (\ref{m4}), (\ref{m13}) and (\ref{m14}).\\
  In the opposite case, upon enlarging $c_{16}$ if necessary we may assume that
  \bas
    M_k \le c_{16} b^k M_{k-1}^{\kappa_k} \qquad \mbox{for all } k\ge 1.
  \eas
  By a straightforward induction, this yields
  \bas
    M_k \le c_{16}^{1+\sum\limits_{j=0}^{k-2} \prod\limits_{i=k-j}^{k} \kappa_i} \cdot
    b^{k+\sum\limits_{j=0}^{k-2} (k-j-1) \cdot \prod\limits_{i=k-j}^{k} \kappa_i} \cdot
    M_0^{\prod\limits_{i=1}^{k} \kappa_i}
  \eas
  for all $k\ge 2$, and hence in view of (\ref{m13}) and (\ref{m5}) we obtain
  \bas
    M_k^\frac{1}{p_k} &\le&
    c_{16}^{\frac{1}{c_1}(\frac{s}{2})^{k}+\frac{1}{c_1} \cdot \sum\limits_{j=0}^{k-2} (\frac{s}{2})^{k-j-1}
    \cdot \prod\limits_{i=k-j}^{k} (1+\eps_i)} \\
    & & \times
    b^{\frac{1}{c_1}k(\frac{s}{2})^{k}+\frac{1}{c_1} \cdot \sum\limits_{j=0}^{k-2} (k-j-1)
        \cdot (\frac{s}{2})^{k-j-1} \cdot \prod\limits_{i=k-j}^{k} (1+\eps_i)} \\
    & & \times
    M_0^{\frac{1}{c_1} \cdot \prod\limits_{i=1}^{k} (1+\eps_i)}
  \eas
  for $k\ge 2$. Since $\ln (1+z) \le z$ for $z\ge 0$, from (\ref{m14}) and the fact that $s<2$ we gain
  \bas
    \ln \Big( \prod_{i=1}^{k} (1+\eps_i) \Big) = \sum_{i=1}^{k} \eps_i \le \frac{c_{18}}{1-\frac{s}{2}},
  \eas
  so that using $\sum\limits_{j=0}^{k-2} (k-j-1) \cdot (\frac{s}{2})^{k-j-1}
    \le \sum\limits_{l=1}^\infty l (\frac{s}{2})^l<\infty$,

  from this we conclude that also in this case $\|\hu(t)\|_{L^\infty(\Omega)}$ is bounded from above by a constant independent
  of $t\in (0,T)$. This clearly proves the lemma.
%
\qed

\vspace*{10mm}
{\bf Acknowledgment.}
  M.~Winkler is grateful for the kind hospitality during his visit at Dong Hua University in March 2010, where
  this work was initiated.\abs
\end{document}